\documentclass[12pt,thmsa]{article}
\usepackage{sw20bams}

\input tcilatex
\begin{document}

\author{Steven Finch, Zai-Qiao Bai and Pascal Sebah}
\title{Typical Dispersion and Generalized Lyapunov Exponents}
\date{March 18, 2008}
\maketitle

\begin{abstract}
Let $f(n)$ denote the number of odd entries in the $n^{\text{th}}$ row of
Pascal's binomial triangle. We study ``average dispersion'' and ``typical
dispersion'' of $f(n)$ -- the latter involves computing a generalized
Lyapunov exponent -- and then turn to numerical analysis of higher
dimensional examples.
\end{abstract}

\footnotetext{
Copyright \copyright\ 2008 by Steven R. Finch. All rights reserved.}Let $D_0$
and $D_1$ denote $m\times m$ nonnegative matrices. Let $z_0$, $z_1$, $\ldots 
$, $z_{k-2}$, $z_{k-1}$ denote a sequence of independent random coin tosses
(heads=1 and tails=0 with equal probability). The Lyapunov exponent
corresponding to random products of $D_0$ and $D_1$ is 
\[
\lambda =\lim_{k\rightarrow \infty }\frac 1k\ln \left\| D_{z_0}D_{z_1}\cdots
D_{z_{k-2}}D_{z_{k-1}}\right\| 
\]
almost surely. We computed $\lambda $ in an earlier paper \cite{FSB} for a
number of examples; our purpose was to quantify the ``typical growth'' of
certain number-theoretic functions. To assess the corresponding ``typical
dispersion'', a second-order extension of $\lambda $ is needed.

The \textbf{generalized Lyapunov exponent} or \textbf{moment Lyapunov
exponent} is 
\[
L(t)=\lim_{k\rightarrow \infty }\frac 1k\ln \left( \limfunc{E}\left( \left\|
D_{z_0}D_{z_1}\cdots D_{z_{k-2}}D_{z_{k-1}}\right\| ^t\right) \right) 
\]
for real $t$. Clearly $L(0)=0$, that is, $e^{L(0)}=1$. Also, $e^{L(1)}$ is
the largest eigenvalue in modulus of $\frac 12(D_0+D_1)$, and $e^{L(2)}$ is
the largest eigenvalue of $\frac 12(D_0\otimes D_0+D_1\otimes D_1)$, where $%
\otimes $ is the direct or Kronecker product of matrices. The latter result
is called the replica trick \cite{CPV}, which can be applied for arbitrary
integer $t\geq 3$ as well.

For simplicity, write $z=z_0z_1\ldots z_{k-2}z_{k-1}$ and $%
D_z=D_{z_0}D_{z_1}\cdots D_{z_{k-2}}D_{z_{k-1}}$. Differentiating $L(t)$
with respect to $t$, we obtain 
\[
L^{\prime }(0)=\lim_{k\rightarrow \infty }\frac 1k\limfunc{E}\left( \ln
\left( \left\| D_z\right\| \right) \right) =\lambda , 
\]
\[
L^{^{\prime \prime }}(0)=\lim_{k\rightarrow \infty }\frac 1k\left[ \limfunc{E%
}\left( \ln \left( \left\| D_z\right\| \right) ^2\right) -\limfunc{E}\left(
\ln \left( \left\| D_z\right\| \right) \right) ^2\right] =\lim_{k\rightarrow
\infty }\frac 1k\limfunc{Var}\left( \ln \left( \left\| D_z\right\| \right)
\right) =\sigma ^2. 
\]
This definition of $\sigma ^2$, unfortunately, is not useful for numerical
calculation.

All $D_0$ matrices exhibited in this paper satisfy $\limfunc{rank}(D_0^q)=1$
for some positive integer $q$. Further, there is an invertible $m\times m$
matrix $Q$ such that $Q^{-1}D_0^qQ$ is the matrix whose $(0,0)^{\text{th}}$
entry is $1$ and all of whose other entries are $0$. Define 
\[
\begin{array}{ccc}
D_0^{\prime }=Q^{-1}D_0Q, &  & D_1^{\prime }=Q^{-1}D_1Q;
\end{array}
\]
by $D_z^{\prime }(0,0)$ is meant the upper left corner entry of $D_z^{\prime
}$. Let $\chi (0^q)$ denote the set of all finite binary words $z$ with no
subwords $0^q$ and with rightmost digit $1$.

Moshe \cite{M6} proved that 
\[
\lambda =\frac 1{2^{q+1}(2^q-1)}\dsum\limits_{z\in \chi (0^q)}\frac
1{2^{\ell (z)}}\ln \left| D_z^{\prime }(0,0)\right| 
\]
where $\ell (z)$ is the length of $z$. We will prove in section [\ref{Prf}]
that 
\[
\sigma ^2=\left( 1+2\frac{2^{2q+1}-(3+q)2^q+1}{2^q-1}\right) \lambda
^2-2\lambda \,\kappa +\mu 
\]
where 
\[
\kappa =\frac 1{2^{q+1}(2^q-1)}\dsum\limits_{z\in \chi (0^q)}\frac{q+\ell (z)%
}{2^{\ell (z)}}\ln \left| D_z^{\prime }(0,0)\right| , 
\]
\[
\mu =\frac 1{2^{q+1}(2^q-1)}\dsum\limits_{z\in \chi (0^q)}\frac 1{2^{\ell
(z)}}\left( \ln \left| D_z^{\prime }(0,0)\right| \right) ^2. 
\]
Summation of the series, coupled with Wynn's $\varepsilon $-process for
accelerating convergence, serves as our primary method for calculating $%
\sigma ^2$. For now, we revisit number-theoretic functions in \cite{FSB} and
compute both ``average dispersion parameters'' $L(2)/\ln (2)$ and ``typical
dispersion parameters'' $\sigma ^2/\ln (2)$.

\subsection{Binomials}

Define $f(n)$ to be the number of odd coefficients in $(1+x)^n$. Let $N$
denote a uniform random integer between $0$ and $n-1$. We have $q=1$, $%
D_0=(1)$, $D_1=(2)$, 
\[
\begin{array}{ccccc}
D_0\otimes D_0=\left( 1\right) , &  & D_1\otimes D_1=\left( 4\right) , &  & 
e^{L(2)}=5/2
\end{array}
\]
hence 
\[
\dfrac{\ln (\limfunc{Var}(f(N)))}{\ln (n)}\rightarrow \dfrac{\ln (5/2)}{\ln
(2)}=1.3219280948873623478703194... 
\]
as $n\rightarrow \infty $. Also $D_0^{\prime }=(1)$, $D_1^{\prime }=(2)$, 
\begin{eqnarray*}
\lambda &=&\frac 14\dsum\limits_{k=0}^\infty \frac 1{2^k}\ln \left|
(D_1^{\prime })^k(0,0)\right| =\frac{\ln (2)}4\dsum\limits_{k=0}^\infty
\frac k{2^k}=\frac{\ln (2)}2 \\
\ &=&0.3465735902799726547086160...,
\end{eqnarray*}
\[
\kappa =\frac 14\dsum\limits_{k=0}^\infty \frac{k+1}{2^k}\ln \left|
(D_1^{\prime })^k(0,0)\right| =\frac{\ln (2)}4\dsum\limits_{k=0}^\infty 
\frac{k(k+1)}{2^k}=2\ln (2), 
\]
\[
\mu =\frac 14\dsum\limits_{k=0}^\infty \frac 1{2^k}\left( \ln \left|
(D_1^{\prime })^k(0,0)\right| \right) ^2=\frac{\ln (2)^2}4\dsum%
\limits_{k=0}^\infty \frac{k^2}{2^k}=\frac 32\ln (2)^2, 
\]
\begin{eqnarray*}
\sigma ^2 &=&3\lambda ^2-2\lambda \,\kappa +\mu =\left( \frac 34-2+\frac
32\right) \ln (2)^2=\frac{\ln (2)^2}4 \\
&=&0.1201132534795503561667756...
\end{eqnarray*}
hence 
\[
\frac{\limfunc{Var}(\ln (f(N)))}{\ln (n)}\rightarrow \frac{\ln (2)}%
4=0.1732867951399863273543080... 
\]
as $n\rightarrow \infty $, confirming a result of Kirschenhofer \cite{Ki}.
The typical parameter $0.173...$ is considerably smaller than the average
parameter $1.321...$ because outlying values (which occur rarely) have been
damped by the logarithm.

\subsection{Trinomials I}

Define $g(n)$ to be the number of odd coefficients in $(1+x+x^2)^n$. We have 
$q=1$, 
\[
\begin{array}{ccc}
D_0=\left( 
\begin{array}{cc}
1 & 2 \\ 
0 & 0
\end{array}
\right) , &  & D_1=\left( 
\begin{array}{cc}
1 & 2 \\ 
1 & 0
\end{array}
\right) ,
\end{array}
\]
\[
\begin{array}{ccc}
D_0\otimes D_0=\left( 
\begin{array}{cccc}
1 & 2 & 2 & 4 \\ 
0 & 0 & 0 & 0 \\ 
0 & 0 & 0 & 0 \\ 
0 & 0 & 0 & 0
\end{array}
\right) , &  & D_1\otimes D_1=\left( 
\begin{array}{cccc}
1 & 2 & 2 & 4 \\ 
1 & 0 & 2 & 0 \\ 
1 & 2 & 0 & 0 \\ 
1 & 0 & 0 & 0
\end{array}
\right)
\end{array}
\]
hence 
\[
\dfrac{\ln (\limfunc{Var}(g(N)))}{\ln (n)}\rightarrow \dfrac{\ln (\xi )}{\ln
(2)}=1.4924205743549514375202537... 
\]
where $\xi =e^{L(2)}=2.813...$ has minimal polynomial $\xi ^3-2\xi ^2-3\xi
+2 $. Also \cite{M5} 
\[
\begin{array}{ccc}
D_0^{\prime }=\left( 
\begin{array}{cc}
1 & 0 \\ 
0 & 0
\end{array}
\right) , &  & D_1^{\prime }=\left( 
\begin{array}{cc}
3 & -4 \\ 
1 & -2
\end{array}
\right) ,
\end{array}
\]
\[
\lambda =\frac 14\dsum\limits_{k=0}^\infty \frac 1{2^k}\ln \left( \frac{%
2^{k+2}-(-1)^k}3\right) =0.4299474333424527201146970..., 
\]
\[
\begin{array}{ccc}
\kappa =\dfrac 14\dsum\limits_{k=0}^\infty \dfrac{k+1}{2^k}\ln \left( \dfrac{%
2^{k+2}-(-1)^k}3\right) , &  & \mu =\dfrac 14\dsum\limits_{k=0}^\infty
\dfrac 1{2^k}\ln \left( \dfrac{2^{k+2}-(-1)^k}3\right) ^2,
\end{array}
\]
\[
\sigma ^2=3\lambda ^2-2\lambda \,\kappa +\mu =0.1211367118847285164803949... 
\]
and 
\[
\frac{\limfunc{Var}(\ln (g(N)))}{\ln (n)}\rightarrow \frac{\sigma ^2}{\ln (2)%
}=0.1747633335056929866262498... 
\]
as $n\rightarrow \infty $.

\subsection{Quadrinomials}

Define $g_3(n)$ to be the number of odd coefficients in $(1+x+x^2+x^3)^n$.
This extends our earlier definitions $f=g_1$ and $g=g_2$. We have $q=2$, 
\[
\begin{array}{ccc}
D_0=\left( 
\begin{array}{ccc}
1 & 2 & 0 \\ 
0 & 0 & 1 \\ 
0 & 0 & 0
\end{array}
\right) , &  & D_1=\left( 
\begin{array}{ccc}
0 & 0 & 0 \\ 
2 & 0 & 0 \\ 
0 & 1 & 2
\end{array}
\right)
\end{array}
\]
hence $e^{L(2)}=5/2$. Also 
\[
\begin{array}{ccc}
D_0^{\prime }=\left( 
\begin{array}{ccc}
1 & 0 & 0 \\ 
0 & 0 & 0 \\ 
0 & 1 & 0
\end{array}
\right) , &  & D_1^{\prime }=\left( 
\begin{array}{ccc}
4 & -4 & -6 \\ 
0 & 2 & 1 \\ 
2 & -4 & -4
\end{array}
\right)
\end{array}
\]
hence \cite{FSB} 
\[
\lambda =\frac 1{24}\dsum\limits_{z\in \chi (00)}\frac 1{2^{\ell (z)}}\ln
\left| D_z^{\prime }(0,0)\right| =\frac{\ln (2)}2, 
\]
\[
\begin{array}{ccc}
\kappa =\dfrac 1{24}\dsum\limits_{z\in \chi (00)}\dfrac{2+\ell (z)}{2^{\ell
(z)}}\ln \left| D_z^{\prime }(0,0)\right| , &  & \mu =\dfrac
1{24}\dsum\limits_{z\in \chi (00)}\dfrac 1{2^{\ell (z)}}\left( \ln \left|
D_z^{\prime }(0,0)\right| \right) ^2,
\end{array}
\]
\[
\sigma ^2=\tfrac{29}3\lambda ^2-2\lambda \,\kappa +\mu =0.12011325..., 
\]
and $\sigma ^2/\ln (2)=0.17328679....$ We conjecture that $\sigma ^2/\ln (2)$
equals $\ln (2)/4$ and prove this to be true in section [\ref{Prf2}].

\subsection{Trinomials II}

Define $h_3(n)$ to be the number of odd coefficients in $(1+x+x^3)^n$. This
extends our earlier definition $g=h_2$. We have $q=2$, 
\[
\begin{array}{ccc}
D_0=\left( 
\begin{array}{cccc}
1 & 2 & 1 & 0 \\ 
0 & 0 & 1 & 1 \\ 
0 & 0 & 0 & 0 \\ 
0 & 0 & 0 & 0
\end{array}
\right) , &  & D_1=\left( 
\begin{array}{cccc}
1 & 1 & 1 & 0 \\ 
1 & 0 & 0 & 1 \\ 
0 & 1 & 0 & 0 \\ 
0 & 0 & 1 & 1
\end{array}
\right)
\end{array}
\]
hence 
\[
\dfrac{\ln (\limfunc{Var}(h_3(N)))}{\ln (n)}\rightarrow \dfrac{\ln (\xi )}{%
\ln (2)}=1.5459492845008943975543991... 
\]
where $\xi =e^{L(2)}=2.919...$ has minimal polynomial 
\[
16\xi ^{10}-40\xi ^9-36\xi ^8+22\xi ^7+76\xi ^6+7\xi ^5-19\xi ^4-19\xi
^3+2\xi +1. 
\]
Also 
\[
\begin{array}{ccc}
D_0^{\prime }=\left( 
\begin{array}{cccc}
1 & 0 & 0 & 0 \\ 
0 & 0 & 0 & 0 \\ 
0 & 1 & 0 & 0 \\ 
0 & 0 & 0 & 0
\end{array}
\right) , &  & D_1^{\prime }=\left( 
\begin{array}{cccc}
3 & -6 & -2 & 4 \\ 
0 & 1 & 1 & 0 \\ 
1 & -3 & -2 & 2 \\ 
0 & 1 & 0 & 0
\end{array}
\right)
\end{array}
\]
hence $\lambda =0.45454538229305...$, $\sigma ^2=0.12497319...$ and 
\[
\frac{\limfunc{Var}(\ln (h_3(N)))}{\ln (n)}\rightarrow \frac{\sigma ^2}{\ln
(2)}=0.18029820... 
\]
as $n\rightarrow \infty $.

\subsection{Quintinomials}

Define $g_4(n)$ to be the number of odd coefficients in $(1+x+x^2+x^3+x^4)^n$%
. We have $q=2$, 
\[
\begin{array}{ccc}
D_0=\left( 
\begin{array}{cccc}
1 & 1 & 2 & 0 \\ 
0 & 0 & 0 & 0 \\ 
0 & 1 & 0 & 2 \\ 
0 & 0 & 0 & 0
\end{array}
\right) , &  & D_1=\left( 
\begin{array}{cccc}
0 & 1 & 2 & 0 \\ 
1 & 0 & 0 & 0 \\ 
1 & 0 & 0 & 2 \\ 
0 & 1 & 0 & 0
\end{array}
\right)
\end{array}
\]
hence 
\[
\dfrac{\ln (\limfunc{Var}(g_4(N)))}{\ln (n)}\rightarrow \dfrac{\ln (\xi )}{%
\ln (2)}=1.6534827473445406557431504... 
\]
where $\xi =e^{L(2)}=3.145...$ has minimal polynomial 
\[
4\xi ^{10}-8\xi ^9-21\xi ^8+14\xi ^7-28\xi ^6+126\xi ^5+65\xi ^4+68\xi
^3+48\xi ^2-56\xi -32. 
\]
Also 
\[
\begin{array}{ccc}
D_0^{\prime }=\left( 
\begin{array}{cccc}
1 & 0 & 0 & 0 \\ 
0 & 0 & 0 & 0 \\ 
0 & 1 & 0 & 0 \\ 
0 & 0 & 0 & 0
\end{array}
\right) , &  & D_1^{\prime }=\left( 
\begin{array}{cccc}
5 & -10 & -8 & 4 \\ 
1 & -1 & -2 & -2 \\ 
1 & -3 & -2 & 4 \\ 
0 & 1 & 0 & -2
\end{array}
\right)
\end{array}
\]
hence $\lambda =0.504253705692...$, $\sigma ^2=0.11406217...$ and 
\[
\frac{\limfunc{Var}(\ln (g_4(N)))}{\ln (n)}\rightarrow \frac{\sigma ^2}{\ln
(2)}=0.16455692... 
\]
as $n\rightarrow \infty $.

\subsection{Trinomials III}

Define $h_4(n)$ to be the number of odd coefficients in $(1+x+x^4)^n$. We
have $q=2$, 
\[
\begin{array}{ccc}
D_0=\left( 
\begin{array}{cccccccc}
\text{{\scriptsize 1}} & \text{{\scriptsize 0}} & \text{{\scriptsize 2}} & 
\text{{\scriptsize 0}} & \text{{\scriptsize 1}} & \text{{\scriptsize 2}} & 
\text{{\scriptsize 1}} & \text{{\scriptsize 1}} \\ 
\text{{\scriptsize 0}} & \text{{\scriptsize 0}} & \text{{\scriptsize 0}} & 
\text{{\scriptsize 0}} & \text{{\scriptsize 0}} & \text{{\scriptsize 0}} & 
\text{{\scriptsize 0}} & \text{{\scriptsize 0}} \\ 
\text{{\scriptsize 0}} & \text{{\scriptsize 1}} & \text{{\scriptsize 0}} & 
\text{{\scriptsize 2}} & \text{{\scriptsize 1}} & \text{{\scriptsize 0}} & 
\text{{\scriptsize 1}} & \text{{\scriptsize 1}} \\ 
\text{{\scriptsize 0}} & \text{{\scriptsize 0}} & \text{{\scriptsize 0}} & 
\text{{\scriptsize 0}} & \text{{\scriptsize 0}} & \text{{\scriptsize 0}} & 
\text{{\scriptsize 0}} & \text{{\scriptsize 0}} \\ 
\text{{\scriptsize 0}} & \text{{\scriptsize 0}} & \text{{\scriptsize 0}} & 
\text{{\scriptsize 0}} & \text{{\scriptsize 0}} & \text{{\scriptsize 0}} & 
\text{{\scriptsize 0}} & \text{{\scriptsize 0}} \\ 
\text{{\scriptsize 0}} & \text{{\scriptsize 0}} & \text{{\scriptsize 0}} & 
\text{{\scriptsize 0}} & \text{{\scriptsize 0}} & \text{{\scriptsize 0}} & 
\text{{\scriptsize 0}} & \text{{\scriptsize 0}} \\ 
\text{{\scriptsize 0}} & \text{{\scriptsize 0}} & \text{{\scriptsize 0}} & 
\text{{\scriptsize 0}} & \text{{\scriptsize 0}} & \text{{\scriptsize 0}} & 
\text{{\scriptsize 0}} & \text{{\scriptsize 0}} \\ 
\text{{\scriptsize 0}} & \text{{\scriptsize 0}} & \text{{\scriptsize 0}} & 
\text{{\scriptsize 0}} & \text{{\scriptsize 0}} & \text{{\scriptsize 0}} & 
\text{{\scriptsize 0}} & \text{{\scriptsize 0}}
\end{array}
\right) , &  & D_1=\left( 
\begin{array}{cccccccc}
\text{{\scriptsize 1}} & \text{{\scriptsize 0}} & \text{{\scriptsize 1}} & 
\text{{\scriptsize 0}} & \text{{\scriptsize 0}} & \text{{\scriptsize 1}} & 
\text{{\scriptsize 0}} & \text{{\scriptsize 0}} \\ 
\text{{\scriptsize 1}} & \text{{\scriptsize 0}} & \text{{\scriptsize 0}} & 
\text{{\scriptsize 0}} & \text{{\scriptsize 0}} & \text{{\scriptsize 0}} & 
\text{{\scriptsize 0}} & \text{{\scriptsize 0}} \\ 
\text{{\scriptsize 0}} & \text{{\scriptsize 1}} & \text{{\scriptsize 0}} & 
\text{{\scriptsize 1}} & \text{{\scriptsize 1}} & \text{{\scriptsize 0}} & 
\text{{\scriptsize 0}} & \text{{\scriptsize 1}} \\ 
\text{{\scriptsize 0}} & \text{{\scriptsize 1}} & \text{{\scriptsize 0}} & 
\text{{\scriptsize 0}} & \text{{\scriptsize 0}} & \text{{\scriptsize 0}} & 
\text{{\scriptsize 0}} & \text{{\scriptsize 0}} \\ 
\text{{\scriptsize 0}} & \text{{\scriptsize 0}} & \text{{\scriptsize 1}} & 
\text{{\scriptsize 0}} & \text{{\scriptsize 0}} & \text{{\scriptsize 0}} & 
\text{{\scriptsize 0}} & \text{{\scriptsize 1}} \\ 
\text{{\scriptsize 0}} & \text{{\scriptsize 0}} & \text{{\scriptsize 0}} & 
\text{{\scriptsize 1}} & \text{{\scriptsize 0}} & \text{{\scriptsize 0}} & 
\text{{\scriptsize 1}} & \text{{\scriptsize 0}} \\ 
\text{{\scriptsize 0}} & \text{{\scriptsize 0}} & \text{{\scriptsize 0}} & 
\text{{\scriptsize 0}} & \text{{\scriptsize 1}} & \text{{\scriptsize 0}} & 
\text{{\scriptsize 0}} & \text{{\scriptsize 0}} \\ 
\text{{\scriptsize 0}} & \text{{\scriptsize 0}} & \text{{\scriptsize 0}} & 
\text{{\scriptsize 0}} & \text{{\scriptsize 0}} & \text{{\scriptsize 1}} & 
\text{{\scriptsize 1}} & \text{{\scriptsize 0}}
\end{array}
\right)
\end{array}
\]
hence 
\[
\dfrac{\ln (\limfunc{Var}(h_4(N)))}{\ln (n)}\rightarrow \dfrac{\ln (\xi )}{%
\ln (2)}=1.5707744868006419128591802... 
\]
where $\xi =e^{L(2)}=2.970...$ has minimal polynomial 
\begin{eqnarray*}
&&32\xi ^{13}-80\xi ^{12}-8\xi ^{11}-60\xi ^{10}-232\xi ^9+240\xi ^8 \\
&&+44\xi ^7+9\xi ^6+11\xi ^5-54\xi ^4-4\xi ^3+3\xi ^2+\xi +2.
\end{eqnarray*}
Also 
\[
\begin{array}{ccc}
D_0^{\prime }=\left( 
\begin{array}{cccccccc}
\text{{\scriptsize 1}} & \text{{\scriptsize 0}} & \text{{\scriptsize 0}} & 
\text{{\scriptsize 0}} & \text{{\scriptsize 0}} & \text{{\scriptsize 0}} & 
\text{{\scriptsize 0}} & \text{{\scriptsize 0}} \\ 
\text{{\scriptsize 0}} & \text{{\scriptsize 0}} & \text{{\scriptsize 0}} & 
\text{{\scriptsize 0}} & \text{{\scriptsize 0}} & \text{{\scriptsize 0}} & 
\text{{\scriptsize 0}} & \text{{\scriptsize 0}} \\ 
\text{{\scriptsize 0}} & \text{{\scriptsize 1}} & \text{{\scriptsize 0}} & 
\text{{\scriptsize 0}} & \text{{\scriptsize 0}} & \text{{\scriptsize 0}} & 
\text{{\scriptsize 0}} & \text{{\scriptsize 0}} \\ 
\text{{\scriptsize 0}} & \text{{\scriptsize 0}} & \text{{\scriptsize 0}} & 
\text{{\scriptsize 0}} & \text{{\scriptsize 0}} & \text{{\scriptsize 0}} & 
\text{{\scriptsize 0}} & \text{{\scriptsize 0}} \\ 
\text{{\scriptsize 0}} & \text{{\scriptsize 0}} & \text{{\scriptsize 0}} & 
\text{{\scriptsize 0}} & \text{{\scriptsize 0}} & \text{{\scriptsize 0}} & 
\text{{\scriptsize 0}} & \text{{\scriptsize 0}} \\ 
\text{{\scriptsize 0}} & \text{{\scriptsize 0}} & \text{{\scriptsize 0}} & 
\text{{\scriptsize 0}} & \text{{\scriptsize 0}} & \text{{\scriptsize 0}} & 
\text{{\scriptsize 0}} & \text{{\scriptsize 0}} \\ 
\text{{\scriptsize 0}} & \text{{\scriptsize 0}} & \text{{\scriptsize 0}} & 
\text{{\scriptsize 0}} & \text{{\scriptsize 0}} & \text{{\scriptsize 0}} & 
\text{{\scriptsize 0}} & \text{{\scriptsize 0}} \\ 
\text{{\scriptsize 0}} & \text{{\scriptsize 0}} & \text{{\scriptsize 0}} & 
\text{{\scriptsize 0}} & \text{{\scriptsize 0}} & \text{{\scriptsize 0}} & 
\text{{\scriptsize 0}} & \text{{\scriptsize 0}}
\end{array}
\right) , &  & D_1^{\prime }=\left( 
\begin{array}{cccccccc}
\text{{\scriptsize 3}} & \text{{\scriptsize 0}} & \text{{\scriptsize -2}} & 
\text{{\scriptsize -8}} & \text{{\scriptsize -4}} & \text{{\scriptsize -2}}
& \text{{\scriptsize -4}} & \text{{\scriptsize -4}} \\ 
\text{{\scriptsize 1}} & \text{{\scriptsize 0}} & \text{{\scriptsize -1}} & 
\text{{\scriptsize -4}} & \text{{\scriptsize -2}} & \text{{\scriptsize -1}}
& \text{{\scriptsize -2}} & \text{{\scriptsize -2}} \\ 
\text{{\scriptsize 0}} & \text{{\scriptsize 1}} & \text{{\scriptsize 0}} & 
\text{{\scriptsize -1}} & \text{{\scriptsize 0}} & \text{{\scriptsize 0}} & 
\text{{\scriptsize -1}} & \text{{\scriptsize 0}} \\ 
\text{{\scriptsize 0}} & \text{{\scriptsize 1}} & \text{{\scriptsize 0}} & 
\text{{\scriptsize -2}} & \text{{\scriptsize -1}} & \text{{\scriptsize 0}} & 
\text{{\scriptsize -1}} & \text{{\scriptsize -1}} \\ 
\text{{\scriptsize 0}} & \text{{\scriptsize 0}} & \text{{\scriptsize 1}} & 
\text{{\scriptsize 0}} & \text{{\scriptsize 0}} & \text{{\scriptsize 0}} & 
\text{{\scriptsize 0}} & \text{{\scriptsize 1}} \\ 
\text{{\scriptsize 0}} & \text{{\scriptsize 0}} & \text{{\scriptsize 0}} & 
\text{{\scriptsize 1}} & \text{{\scriptsize 0}} & \text{{\scriptsize 0}} & 
\text{{\scriptsize 1}} & \text{{\scriptsize 0}} \\ 
\text{{\scriptsize 0}} & \text{{\scriptsize 0}} & \text{{\scriptsize 0}} & 
\text{{\scriptsize 0}} & \text{{\scriptsize 1}} & \text{{\scriptsize 0}} & 
\text{{\scriptsize 0}} & \text{{\scriptsize 0}} \\ 
\text{{\scriptsize 0}} & \text{{\scriptsize 0}} & \text{{\scriptsize 0}} & 
\text{{\scriptsize 0}} & \text{{\scriptsize 0}} & \text{{\scriptsize 1}} & 
\text{{\scriptsize 1}} & \text{{\scriptsize 0}}
\end{array}
\right)
\end{array}
\]
hence $\lambda =0.45759385431410...$, $\sigma ^2=0.13055386...$ and 
\[
\frac{\limfunc{Var}(\ln (h_4(N)))}{\ln (n)}\rightarrow \frac{\sigma ^2}{\ln
(2)}=0.18834940... 
\]
as $n\rightarrow \infty $.

\subsection{Sextinomials}

Define $g_5(n)$ to be the number of odd coefficients in $(1+x+\ldots
+x^4+x^5)^n$. We have $q=3$, 
\[
\begin{array}{ccc}
D_0=\left( 
\begin{array}{cccccc}
1 & 1 & 2 & 2 & 0 & 0 \\ 
0 & 0 & 0 & 0 & 0 & 0 \\ 
0 & 1 & 0 & 0 & 1 & 1 \\ 
0 & 0 & 0 & 0 & 0 & 0 \\ 
0 & 0 & 0 & 0 & 0 & 0 \\ 
0 & 0 & 0 & 0 & 1 & 0
\end{array}
\right) , &  & D_1=\left( 
\begin{array}{cccccc}
0 & 0 & 0 & 0 & 0 & 0 \\ 
2 & 2 & 0 & 0 & 0 & 0 \\ 
0 & 0 & 0 & 0 & 0 & 0 \\ 
0 & 0 & 1 & 1 & 2 & 2 \\ 
0 & 0 & 0 & 1 & 0 & 0 \\ 
0 & 0 & 0 & 0 & 0 & 0
\end{array}
\right)
\end{array}
\]
hence 
\[
\dfrac{\ln (\limfunc{Var}(g_5(N)))}{\ln (n)}\rightarrow \dfrac{\ln (\xi )}{%
\ln (2)}=1.6903750759639444915537652... 
\]
where $\xi =e^{L(2)}=3.227...$ has minimal polynomial 
\begin{eqnarray*}
&&128\xi ^{11}-640\xi ^{10}+416\xi ^9+1008\xi ^8+416\xi ^7-28\xi ^6 \\
&&-3112\xi ^5-2572\xi ^4+346\xi ^3+1887\xi ^2+511\xi +144.
\end{eqnarray*}
Also 
\[
\begin{array}{ccc}
D_0^{\prime }=\left( 
\begin{array}{cccccc}
1 & 0 & 0 & 0 & 0 & 0 \\ 
0 & 0 & 0 & 0 & 0 & 0 \\ 
0 & 1 & 0 & 0 & 0 & 0 \\ 
0 & 0 & 1 & 0 & 0 & 0 \\ 
0 & 0 & 0 & 0 & 0 & 0 \\ 
0 & 0 & 0 & 0 & 0 & 0
\end{array}
\right) , &  & D_1^{\prime }=\left( 
\begin{array}{cccccc}
6 & -8 & -8 & -10 & -4 & -6 \\ 
0 & 0 & 0 & 0 & 0 & 1 \\ 
2 & -4 & -4 & -4 & 0 & -3 \\ 
0 & 0 & 0 & 0 & 0 & 0 \\ 
2 & -4 & -4 & -4 & 0 & -3 \\ 
0 & 2 & 2 & 1 & -2 & 1
\end{array}
\right)
\end{array}
\]
hence 
\[
\lambda =\frac 1{112}\dsum\limits_{z\in \chi (000)}\frac 1{2^{\ell (z)}}\ln
\left| D_z^{\prime }(0,0)\right| =0.5344481528..., 
\]
\[
\begin{array}{ccc}
\kappa =\dfrac 1{112}\dsum\limits_{z\in \chi (000)}\dfrac{3+\ell (z)}{%
2^{\ell (z)}}\ln \left| D_z^{\prime }(0,0)\right| , &  & \mu =\dfrac
1{112}\dsum\limits_{z\in \chi (000)}\dfrac 1{2^{\ell (z)}}\left( \ln \left|
D_z^{\prime }(0,0)\right| \right) ^2,
\end{array}
\]
\[
\sigma ^2=\tfrac{169}7\lambda ^2-2\lambda \,\kappa +\mu =0.0965... 
\]
and 
\[
\frac{\limfunc{Var}(\ln (g_5(N)))}{\ln (n)}\rightarrow \frac{\sigma ^2}{\ln
(2)}=0.1392... 
\]
as $n\rightarrow \infty $.

\subsection{Septinomials}

Define $g_6(n)$ to be the number of odd coefficients in $(1+x+\cdots
+x^5+x^6)^n$. We have $q=3$, 
\[
\begin{array}{ccc}
D_0=\left( 
\begin{array}{cccccc}
1 & 0 & 1 & 2 & 0 & 0 \\ 
0 & 0 & 0 & 0 & 0 & 0 \\ 
0 & 0 & 0 & 0 & 1 & 2 \\ 
0 & 2 & 1 & 0 & 1 & 0 \\ 
0 & 0 & 0 & 0 & 0 & 0 \\ 
0 & 0 & 0 & 0 & 0 & 0
\end{array}
\right) , &  & D_1=\left( 
\begin{array}{cccccc}
0 & 0 & 0 & 2 & 1 & 0 \\ 
1 & 0 & 0 & 0 & 0 & 0 \\ 
1 & 0 & 0 & 0 & 0 & 2 \\ 
0 & 2 & 1 & 0 & 0 & 0 \\ 
0 & 0 & 1 & 0 & 0 & 0 \\ 
0 & 0 & 0 & 0 & 1 & 0
\end{array}
\right)
\end{array}
\]
hence 
\[
\dfrac{\ln (\limfunc{Var}(g_6(N)))}{\ln (n)}\rightarrow \dfrac{\ln (\xi )}{%
\ln (2)}=1.7258729504941114967801068... 
\]
where $\xi =e^{L(2)}=3.307...$ has minimal polynomial 
\begin{eqnarray*}
&&8\xi ^{21}-4\xi ^{20}-18\xi ^{19}-335\xi ^{18}+34\xi ^{17}+474\xi
^{16}+4072\xi ^{15}+302\xi ^{14} \\
&&-3119\xi ^{13}-16848\xi ^{12}-1056\xi ^{11}+7321\xi ^{10}+29681\xi
^9+910\xi ^8 \\
&&-6690\xi ^7-22628\xi ^6-152\xi ^5+1936\xi ^4+6112\xi ^3-128\xi -512.
\end{eqnarray*}
Also 
\[
\begin{array}{ccc}
D_0^{\prime }=\left( 
\begin{array}{cccccc}
1 & 0 & 0 & 0 & 0 & 0 \\ 
0 & 0 & 0 & 0 & 0 & 0 \\ 
0 & 1 & 0 & 0 & 0 & 0 \\ 
0 & 0 & 1 & 0 & 0 & 0 \\ 
0 & 0 & 0 & 0 & 0 & 0 \\ 
0 & 0 & 0 & 0 & 0 & 0
\end{array}
\right) , &  & D_1^{\prime }=\left( 
\begin{array}{cccccc}
7 & -36 & -28 & -24 & 4 & -4 \\ 
0 & 0 & 1 & 0 & \frac 12 & -2 \\ 
\frac 32 & -8 & -8 & -6 & \frac 12 & 1 \\ 
0 & 0 & 1 & 0 & -\frac 12 & 1 \\ 
2 & -12 & -10 & -8 & 1 & 0 \\ 
1 & -6 & -6 & -4 & 1 & 0
\end{array}
\right)
\end{array}
\]
hence $\lambda =0.53765282...$, $\sigma ^2=0.1082...$ and 
\[
\frac{\limfunc{Var}(\ln (g_6(N)))}{\ln (n)}\rightarrow \frac{\sigma ^2}{\ln
(2)}=0.1561... 
\]
as $n\rightarrow \infty $.

\subsection{\label{Prf}Proof of Formula for $\sigma ^2$}

Fix a nonnegative integer $k$. In the definitions of $\lambda $ and $\sigma
^2$, we assumed that each binary word $z$ of length $k$ occurs with
probability $2^{-k}$. This assumption is not necessary: let $p_k(z)$ denote
the (non-uniform) probability associated with $z$. Let 
\[
\begin{array}{ccc}
c_k(t)=\dsum\limits_{\ell (z)=k}p_k(z)\left\| D_z\right\| ^t, &  & t\text{
real}
\end{array}
\]
then the radius of convergence of 
\[
\begin{array}{ccc}
\Omega (s,t)=\dsum\limits_zs^{\ell (z)}p_{\ell (z)}(z)\left\| D_z\right\|
^t=\dsum\limits_{k=0}^\infty c_k(t)s^k, &  & s\text{ complex}
\end{array}
\]
is 
\[
\lim_{k\rightarrow \infty }c_k^{-1/k}=\lim_{k\rightarrow \infty }\left( 
\limfunc{E}\left( \left\| D_z\right\| ^t\right) \right) ^{-1/k}=e^{-L(t)}. 
\]
Hence $\Omega (s,t)$ converges absolutely if $|s|<e^{-L(t)}$ and diverges if 
$|s|>e^{-L(t)}$. For convenience, we will write $p(z)$ instead of $p_{\ell
(z)}(z)$ from now on.

We also postulated that there exists a positive integer $q$ such that $%
\limfunc{rank}(D_0^q)=1$. This postulate can be weakened to the following:
there exists a binary word $z^{*}$ for which $\limfunc{rank}(D_{z^{*}})=1$.
Any product of matrices involving $D_{z^{*}}$ is also rank $1$ since, if $%
D_{z^{*}}=\alpha \,\beta ^T$, then 
\[
D_xD_{z^{*}}D_y=(D_x\alpha )(\beta ^TD_y). 
\]
Let $\chi (z^{*})$ denote the set of all finite binary words $w$ such that $%
z^{*}$ appears in the word $w\,z^{*}$ only at the end. Hence, given an
arbitrary binary word $z$, it follows that either 
\[
z=w^{(0)} 
\]
for some $w^{(0)}\in \chi (z^{*})$ or 
\[
z=(w^{(\infty )}z^{*})(w^{(j)}z^{*})(w^{(j-1)}z^{*})\cdots
(w^{(2)}z^{*})(w^{(1)}z^{*})w^{(0)} 
\]
uniquely, where $j\geq 0$ and $w^{(i)}\in \chi (z^{*})$ for all $i$. In the
latter case, clearly 
\[
\left\| D_z\right\| =|\beta ^TD_{w^{(1)}}\alpha |\cdot |\beta
^TD_{w^{(2)}}\alpha |\cdots |\beta ^TD_{w^{(j-1)}}\alpha |\cdot |\beta
^TD_{w^{(j)}}\alpha |\cdot \left\| (D_{w^{(\infty )}}\alpha )(\beta
^TD_{w^{(0)}})\right\| . 
\]
Consider the words $w^{(0)}$, $w^{(\infty )}$ as fixed and the index $j$ as
increasing; the final factor is thus immaterial. It follows that the
convergence behavior of $\Omega (s,t)$ is the same as 
\begin{eqnarray*}
&&\ \ \dsum\limits_{w^{(1)}}s^{\ell (w^{(1)}z^{*})}p(w^{(1)}z^{*})|\beta
^TD_{w^{(1)}}\alpha |^t+ \\
&&\ \ \dsum\limits_{w^{(2)}}\dsum\limits_{w^{(1)}}s^{\ell
(w^{(2)}z^{*}w^{(1)}z^{*})}p(w^{(2)}z^{*}w^{(1)}z^{*})|\beta
^TD_{w^{(2)}}\alpha |^t|\beta ^TD_{w^{(1)}}\alpha |^t+ \\
&&\ \
\dsum\limits_{w^{(3)}}\dsum\limits_{w^{(2)}}\dsum\limits_{w^{(1)}}s^{\ell
(w^{(3)}z^{*}w^{(2)}z^{*}w^{(1)}z^{*})}p(w^{(3)}z^{*}w^{(2)}z^{*}w^{(1)}z^{*})|\beta ^TD_{w^{(3)}}\alpha |^t|\beta ^TD_{w^{(2)}}\alpha |^t|\beta ^TD_{w^{(1)}}\alpha |^t+\cdots
\end{eqnarray*}
which equals $F(s,t)/(1-F(s,t)),$ where 
\[
F(s,t)=\dsum\limits_{w\in \chi (z^{*})}s^{\ell (w\,z^{*})}p(w\,z^{*})|\beta
^TD_w\alpha |^t. 
\]
Therefore the smallest zero in modulus of $1-F(s,t)$ determines the radius
of convergence of $\Omega (s,t)$ and the generalized Lyapunov exponent
satisfies 
\[
\begin{array}{ccccc}
L(t)=-\ln (s(t)) &  & \text{where} &  & F(s(t),t)=1.
\end{array}
\]

From the formula for $L(t)$ in terms of $s(t)$, we deduce that 
\[
s(0)=1 
\]
since $L(0)=0$, thus $F(1,0)=1$. By the Implicit Function Theorem, 
\[
s^{\prime }(0)=\left. -\frac{\partial F(s,t)/\partial t}{\partial
F(s,t)/\partial s}\right| \Sb s=1  \\ t=0  \endSb =-\frac{F_t(1,0)}{F_s(1,0)}
\]
thus 
\[
\lambda =L^{\prime }(0)=-\frac{s^{\prime }(0)}{s(0)}=\frac{F_t(1,0)}{F_s(1,0)%
}. 
\]
By the Chain Rule, 
\begin{eqnarray*}
s^{\prime \prime }(0) &=&\left. -\frac d{dt}\frac{F_t(s,t)}{F_s(s,t)}\right| 
\Sb s=1  \\ t=0  \endSb  \\
\ &=&\left. -\frac \partial {\partial s}\frac{F_t(s,t)}{F_s(s,t)}\right| \Sb %
s=1  \\ t=0  \endSb \cdot s^{\prime }(0)\left. -\frac \partial {\partial t}%
\frac{F_t(s,t)}{F_s(s,t)}\right| \Sb s=1  \\ t=0  \endSb  \\
\ &=&-\left( \frac{F_{st}}{F_s}-\frac{F_{ss}F_t}{F_s^2}\right) \left( -\frac{%
F_t}{F_s}\right) -\left( \frac{F_{tt}}{F_s}-\frac{F_{st}F_t}{F_s^2}\right)
\end{eqnarray*}
thus 
\begin{eqnarray*}
\sigma ^2 &=&L^{\prime \prime }(0)=-\frac{s(0)s^{\prime \prime
}(0)-s^{\prime }(0)^2}{s(0)^2}=s^{\prime }(0)^2-s^{\prime \prime }(0) \\
\ &=&\frac{(F_{ss}+F_s)F_t^2}{F_s^3}-2\frac{F_{st}F_t}{F_s^2}+\frac{F_{tt}}{%
F_s}.
\end{eqnarray*}
It is easy to show that 
\[
F_t=\dsum\limits_{w\in \chi (z^{*})}p(w\,z^{*})\ln |\beta ^TD_w\alpha |, 
\]
\[
F_{st}=\dsum\limits_{w\in \chi (z^{*})}\ell (w\,z^{*})p(w\,z^{*})\ln |\beta
^TD_w\alpha |, 
\]
\[
F_{tt}=\dsum\limits_{w\in \chi (z^{*})}p(w\,z^{*})\left( \ln |\beta
^TD_w\alpha |\right) ^2. 
\]
In our scenario, $z^{*}=0^q$. The derivatives $F_s$ and $F_{ss}$ can be
written in closed-form because 
\[
\chi (0^q)=\{\emptyset \}\cup \chi _0\cup (\chi _0\times \chi _0)\cup (\chi
_0\times \chi _0\times \chi _0)\cup \cdots 
\]
where $\chi _0=\{1,01,0^21,\ldots ,0^{q-1}1\}$ and the Cartesian product $%
\chi _0\times \chi _0$ is to be interpreted as concatenation: 
\[
\chi _0\times \chi _0= 
\begin{array}{l}
\{11,101,10^21,\ldots ,10^{q-1}1, \\ 
\;\;011,0101,010^21,\ldots ,010^{q-1}1, \\ 
\;\;0^211,0^2101,0^210^21,\ldots ,0^210^{q-1}1, \\ 
\;\;\vdots \\ 
\;\;0^{q-1}11,0^{q-1}101,0^{q-1}10^21,\ldots ,0^{q-1}10^{q-1}1\}
\end{array}
\]
hence 
\begin{eqnarray*}
F(s,0) &=&\dsum\limits_{w\in \chi (0^q)}s^{\ell (w)+q}2^{-\ell (w)-q}=\left(
\frac s2\right) ^q\dsum\limits_{m=0}^\infty \left(
\dsum\limits_{k=1}^q\left( \frac s2\right) ^k\right) ^m \\
\ &=&\left( \frac s2\right) ^q\left( 1-\dsum\limits_{k=1}^q\left( \frac
s2\right) ^k\right) ^{-1}=\left( \frac s2\right) ^q\left( 1-\frac s2\right)
\left( 1-s+\left( \frac s2\right) ^{q+1}\right) ^{-1}.
\end{eqnarray*}
It follows that 
\[
F_s(1,0)=2(2^q-1), 
\]
\[
F_{ss}(1,0)=4(2^{2q+1}-(3+q)2^q+1) 
\]
and therefore 
\begin{eqnarray*}
\sigma ^2 &=&\frac{F_{ss}}{F_s}\lambda ^2+\lambda ^2-2\frac{F_{st}}{F_s}%
\lambda +\frac{F_{tt}}{F_s} \\
&=&\left( 1+2\frac{2^{2q+1}-(3+q)2^q+1}{2^q-1}\right) \lambda ^2-2\lambda 
\frac{F_{st}}{F_s}+\frac{F_{tt}}{F_s}.
\end{eqnarray*}
We conclude the proof by setting $\kappa =F_{st}/F_s$ and $\mu =F_{tt}/F_s$.

\subsection{\label{Prf2}Proof of Conjecture}

Before giving the proof, let us extend the formula for $F(s,t)$ to the case
of three nonnegative square matrices $E_0$, $E_1$, $E_2$ satisfying the
following: there exist two ternary words $z_1^{*}$, $z_2^{*}$ for which $%
\limfunc{rank}(E_{z_1^{*}})=\limfunc{rank}(E_{z_2^{*}})=1$. Write $%
E_{z_1^{*}}=\alpha _1\,\beta _1^T$, $E_{z_2^{*}}=\alpha _2\,\beta _2^T$ and
assume that $z_i^{*}$ is not a subword of $z_j^{*}$ for $i\neq j$. Let $\chi
(z_1^{*},z_2^{*})$ denote the set of all finite ternary words $w$ such that $%
z_i^{*}$ appears in the word $w\,z_i^{*}$ only at the end, $1\leq i\leq 2$.
A unique factorization property holds for all ternary words $z$ as before.
The function $F(s,t)$ becomes a $2\times 2$ matrix with entries 
\[
F_{i,j}(s,t)=\dsum\limits_{w\in \chi (z_1^{*},z_2^{*})}s^{\ell
(w\,z_i^{*})}p(w\,z_i^{*})|\beta _j^TE_w\alpha _i|^t. 
\]
The smallest zero in modulus of $\limfunc{det}(I-F(s,t))$ gives rise to the
generalized Lyapunov exponent $L(t)=-\ln (s(t))$ analogous to before.

Let us return to binary words $z$. We assume without loss of generality that 
$z_0=1$ (for this argument may be repeated with $0$s and $1$s interchanged).
Every word thus looks like 
\[
z=1\,0^{j_0}\,1\,0^{j_1}\,1\,0^{j_2}\,\ldots
\,1\,0^{j_{n-2}}\,1\,0^{j_{n-1}} 
\]
where each $j_i\geq 0$. We introduce a rewording of $D_z$: 
\[
\tilde D_{j_0}\tilde D_{j_1}\tilde D_{j_2}\cdots \tilde D_{j_{n-2}}\tilde
D_{j_{n-1}}=\left( D_1D_0^{j_0}\right) \left( D_1D_0^{j_1}\right) \left(
D_1D_0^{j_2}\right) \cdots \left( D_1D_0^{j_{n-2}}\right) \left(
D_1D_0^{j_{n-1}}\right) 
\]
and note that the weight associated with $\tilde D_{j_i}$ is $(s/2)^{j_i+1}$%
. Also 
\[
\begin{array}{ccc}
\tilde D_0=D_1=\left( 
\begin{array}{ccc}
0 & 0 & 0 \\ 
2 & 0 & 0 \\ 
0 & 1 & 2
\end{array}
\right) , &  & \tilde D_1=D_1D_0=\left( 
\begin{array}{ccc}
0 & 0 & 0 \\ 
2 & 4 & 0 \\ 
0 & 0 & 1
\end{array}
\right) ,
\end{array}
\]
\[
\tilde D_j=D_1D_0^j=\left( 
\begin{array}{ccc}
0 & 0 & 0 \\ 
2 & 4 & 4 \\ 
0 & 0 & 0
\end{array}
\right) 
\]
for all $j\geq 2$. Hence the rewording actually consists of only three
matrices $\tilde D_0$, $\tilde D_1$, $\tilde D_2$ with weights $s/2$, $s^2/4$
and 
\[
\dsum\limits_{m=2}^\infty \left( \frac s2\right) ^{m+1}=\frac{s^3/8}{1-s/2}. 
\]
Further, the initial row of each matrix is zero, thus we may consider only
the lower-right $2\times 2$ submatrix: 
\[
E_0=\left( 
\begin{array}{cc}
0 & 0 \\ 
1 & 2
\end{array}
\right) =\left( 
\begin{array}{c}
0 \\ 
1
\end{array}
\right) \left( 
\begin{array}{cc}
1 & 2
\end{array}
\right) =\alpha _1\beta _1^T, 
\]
\[
E_1=\left( 
\begin{array}{cc}
4 & 0 \\ 
0 & 1
\end{array}
\right) =M, 
\]
\[
E_2=\left( 
\begin{array}{cc}
4 & 4 \\ 
0 & 0
\end{array}
\right) =\left( 
\begin{array}{c}
1 \\ 
0
\end{array}
\right) \left( 
\begin{array}{cc}
4 & 4
\end{array}
\right) =\alpha _2\beta _2^T. 
\]
Let $r_1(s)=s/2$, $r_2(s)=s^2/4$ denote the weights corresponding to $\alpha
_1\beta _1^T$, $\alpha _2\beta _2^T$ and $q(s)=(s^3/8)/(1-s/2)$ denote the
weight corresponding to $M$. Setting $z_1^{*}=0$ and $z_2^{*}=2$, we obtain 
\[
F_{i,j}(s,t)=\dsum\limits_{k=0}^\infty r_i(s)q(s)^k\left| \beta
_j^TM^k\alpha _i\right| ^t. 
\]
As calculated in \cite{FSB}, 
\[
\beta _j^TM^k\alpha _i=\left\{ 
\begin{array}{ccc}
2 &  & \text{if }i=1\text{, }j=1 \\ 
2^2 &  & \text{if }i=1\text{, }j=2 \\ 
2^{2k} &  & \text{if }i=2\text{, }j=1 \\ 
2^{2(k+1)} &  & \text{if }i=2\text{, }j=2
\end{array}
\right. 
\]
hence 
\[
F_{1,1}(s,t)=\dsum\limits_{k=0}^\infty \frac s2\left( \frac{s^2}4\right)
^k2^t=\frac{2^t\zeta }{1-\zeta ^2} 
\]
\[
F_{1,2}(s,t)=\dsum\limits_{k=0}^\infty \frac s2\left( \frac{s^2}4\right)
^k2^{2t}=\frac{2^{2t}\zeta }{1-\zeta ^2} 
\]
\[
F_{2,1}(s,t)=\dsum\limits_{k=0}^\infty \frac{s^3/8}{1-s/2}\left( \frac{s^2}%
4\right) ^k2^{2k\,t}=\frac{\zeta ^3}{(1-\zeta )(1-2^{2t}\zeta ^2)} 
\]
\[
F_{2,2}(s,t)=\dsum\limits_{k=0}^\infty \frac{s^3/8}{1-s/2}\left( \frac{s^2}%
4\right) ^k2^{2(k+1)t}=\frac{2^{2t}\zeta ^3}{(1-\zeta )(1-2^{2t}\zeta ^2)} 
\]
where $\zeta =s/2$. The smallest zero of 
\[
\limfunc{det}(F(s,t)-I)=\frac{\left( 1-(2^t+1)\zeta \right) \left(
1-(1-2^t+2^{2t})\zeta ^2\right) }{(1-\zeta )(1-\zeta ^2)(1-2^{2t}\zeta ^2)} 
\]
is $\zeta (t)=1/(2^t+1)$, hence 
\[
L(t)=-\ln (2\zeta (t))=\ln ((2^t+1)/2). 
\]
This expression for $L(t)$ is trivially the same as the generalized Lyapunov
exponent for binomials, thus $\lambda =\ln (2)/2$ and $\sigma ^2=\ln (2)^2/4$%
.

\subsection{Digital Sums}

Let $\#(n)$ denote the number of $1$s in the binary expansion of $n$. We
know that $\#(n)=\ln (f(n))/\ln (2)$ and $\#(3n)=\ln (g_3(n))/\ln (2)$,
therefore $\#(N)$ and $\#(3N)$ are identically distributed in the sense that 
\[
\limfunc{P}\left( \frac{\#(N)-\tfrac{\ln (n)}{2\ln (2)}}{\sqrt{\tfrac{\ln (n)%
}{4\ln (2)}}}\leq x\right) \rightarrow \frac 1{\sqrt{2\pi }%
}\dint\limits_{-\infty }^x\exp \left( -\frac{t^2}2\right) dt\leftarrow 
\limfunc{P}\left( \frac{\#(3N)-\tfrac{\ln (n)}{2\ln (2)}}{\sqrt{\tfrac{\ln
(n)}{4\ln (2)}}}\leq x\right) 
\]
as $n\rightarrow \infty $. Moreover, $\#(aN+b)$ are identically distributed, 
$0\leq b<a$. This answers a question raised in \cite{FSB} and we refer
interested readers to details in \cite{S1, S2}.

\subsection{Closing Words}

It is natural to seek a Central Limit Theorem for functions examined in this
paper, for example, $\#(n)$ or $f(n)$. This seems to be an open problem, but
we indicate a possible direction for solution. Trollope \cite{Tr} \& Delange 
\cite{De} proved that 
\[
\frac 1n\dsum\limits_{k=0}^{n-1}\#(k)-\dfrac 1{2\ln (2)}\ln (n)=\Phi \left( 
\frac{\ln (n)}{\ln (2)}\right) 
\]
\textit{exactly}, where $\Phi (x)$ is a certain continuous
nowhere-differentiable function of period 1, 
\[
-0.2075187496394219092731305...=\frac{\ln (3)}{2\ln (2)}-1=\inf_x\Phi
(x)<\sup_x\Phi (x)=0, 
\]
and the Fourier coefficients of $\Phi (x)$ are all known. The mean value of $%
\Phi (x)$ is \cite{De, FGKPT} 
\[
\dint\limits_0^1\Phi (x)dx=\frac 1{2\ln (2)}\left( \ln (2\pi )-1\right)
-\frac 34=-0.1455994557083223046583226.... 
\]
In principle, we can numerically compute any percentile of $\Phi (x)$ by
approximating the Lebesgue measure of all $x\in [0,1]$ satisfying $\Phi
(x)<\tau $ for some threshold $\tau $. This ``inversion'' must be done
carefully, however, because of the fractal nature of $\Phi (x)$. A plot of
the limiting density function would be good to see!

Stein \cite{St} \&\ Larcher \cite{La} likewise proved that

\[
\ln \left( \frac 1n\dsum\limits_{k=0}^{n-1}f(k)\right) -\dfrac{\ln (3/2)}{%
\ln (2)}\ln (n)=\ln \left( \Psi \left( \frac{\ln (n)}{\ln (2)}\right)
\right) 
\]
where $\Psi (x)$ is a continuous nowhere-monotonic (but almost-everywhere
differentiable) function of period 1, 
\[
0.8125565590160063876948821...=\inf_x\Psi (x)<\sup_x\Psi (x)=1, 
\]
and the Fourier coefficients of $\Psi (x)$ are all known. No closed-form
expression is evident for the lower bound \cite{Fi, Ws}. The mean value of $%
\Psi (x)$ is 
\[
\dint\limits_0^1\Psi (x)dx=0.8636049963990796049605033... 
\]
and this too is unknown \cite{FGKPT, GH}. Again, a plot of the limiting
density corresponding to $\ln (\Psi (x))$ would be welcome progress.

The function $g(n)$ deserves more attention: Fourier expansions for the
analogs of both $\Phi (x)$ and $\Psi (x)$ are desired (if these exist). We
have not mentioned thus far the functions $u(n)$ associated with Pascal's
rhombus or $v(n)$ associated with ``Fibonacci's rhombus'' \cite{FSB}. For
Stern's sequence, in which $v(n)$ is the number of odd coefficients in 
\[
\begin{array}{ccccc}
p_n(x)=x\,p_{n-1}(x)+p_{n-2}(x), &  & p_1(x)=x, &  & p_0(x)=1
\end{array}
\]
we have 
\[
\lambda =0.396212564297744..., 
\]
\[
\sigma ^2=0.022172945128737... 
\]
and the method for calculating $\sigma ^2$ will be published later \cite{Ba}%
. By contrast, $u(n)$ is the number of odd coefficients in 
\[
\begin{array}{ccccc}
p_n(x)=(1+x+x^2)p_{n-1}(x)+x^2p_{n-2}(x), &  & p_1(x)=1+x+x^2, &  & p_0(x)=1.
\end{array}
\]
and $\lambda =0.57331379313...$, but no precise estimate of $\sigma ^2$ has
yet been found.

\end{document}